\documentclass[a4paper,12pt]{article}
\usepackage{amsmath,amssymb,amsthm,graphicx,floatflt,mathtext,enumitem,pifont,xcolor,epsfig,multirow,pdfpages,bm,multicol}
\usepackage[cp1251]{inputenc}
\usepackage[T2A]{fontenc}
\usepackage[english]{babel}
\usepackage{indentfirst}

\begin{document}
\numberwithin{equation}{section}
\newtheorem{definition}{Definition}[section]
\newtheorem{theorem}{Theorem}[section]
\newtheorem{lemma}{Lemma}[section]
\newtheorem{statement}{Statement}[section]
\newtheorem{proposition}{Proposition}[section]
\newtheorem{remark}{Remark}[section]
\newtheorem{corollary}{Corollary}[section]

\newcommand{\norm}[2]{\left\lVert#1\right\rVert_{#2}} 
\newcommand{\ip}[3]{\left\langle#1,#2\right\rangle_{#3}} 
\newcommand{\abs}[1]{\left\lvert#1\right\rvert} 
\newcommand{\restr}[2]{\left.{#1}\right|_{#2}} 
\newcommand{\intl}{\int\limits_0^\ell}
\newcommand{\dd}[1]{\mathrm{d}#1} 
\newcommand{\dxx}{\dd\zeta} 
\newcommand{\xx}{(\zeta+\ell_0)}
\newcommand{\ra}{\rho A}
\newcommand{\ac}{{\tilde A}} 
\noindent UDC 517.977

\begin{center}\textsc{\large{\textbf{Partial Stabilization of an Orbiting Satellite Model \\
with a Flexible Attachment}}}\end{center}
\setlength{\jot}{10pt}

\begin{center}
{Julia Kalosha\footnote{Institute of Applied Mathematics \& Mechanics, National Academy of Sciences of Ukraine, Sloviansk, Ukraine}, Yevgeniia Yevgenieva$^{1,}$\footnote{Max Planck Institute for Dynamics of Complex Technical Systems, Magdeburg, Germany}, Alexander Zuyev$^{1,2}$}
\end{center}

\begin{abstract}
We consider a mathematical model of an orbiting satellite, comprising a rigid carrier body and a flexible boom, operating under the influence of gravity gradient torque.
This model is represented by a nonlinear control system, which includes ordinary differential equations governing the carrier body's angular velocity and attitude quaternion, coupled with the Euler-Bernoulli equations that describe the vibration of the flexible component.
We propose an explicit feedback design aimed at guaranteeing the partial stability of the closed-loop system in an appropriate Hilbert space.
\end{abstract}

\section{Introduction}

The stabilization problem for satellites with flexible attachments is a crucial challenge in aerospace research and engineering, as these attachments (solar panels, booms, tethers, and antennas) are essential for the satellite's operation but can introduce complex dynamics. Without proper control schemes, these flexible parts can induce vibrations  that degrade the satellite's overall performance and longevity. Effective stabilization ensures accurate orientation and stability, which are vital for tasks like communication, Earth observation, and scientific measurements.

The fundamental applications of Lyapunov's direct method in analyzing the stability of satellites with elastic components are summarized in~\cite{meirovitch1972liapunov,meirovitch1972stability,nabiullin}, while a brief survey of the latest findings in this field is available in~\cite{gordon2023effects}.
The recent study presented in~\cite{gordon2023effects} explores how the model fidelity and the uncertainty of parameters influence the effectiveness of a combined model-based feedback and feedforward control approach for maintaining the orientation of a satellite equipped with flexible appendages.

The present paper focuses on the analytical control design for a realistic model of a flexible satellite model within the framework of partial stability theory~\cite{rumyantsev,zuev2000stabilization}. A related single-axis stabilization problem was solved in~\cite{zuyev2001partial} for rigid satellite models actuated by a pair of thrusters and reaction wheels.
The structure of this paper is arranged in the following manner.
In Section~\ref{sec_model}, we derive a mathematical model of an orbiting satellite with a boom. It is assumed that the satellite consists of a rigid carrier body moving along the Earth's orbit and a boom in the form of an elastic cantilever beam.
We derive the equations describing the vibrations of the boom by applying Hamilton's principle with an account of external forces. A feedback control is proposed in Section~\ref{sec_control} with the use of an energy-based Lyapunov functional.
The corresponding closed-loop system is represented by an abstract differential equation in Section~\ref{sec_abstract}, and the main stability result is formulated in Section~\ref{sec_thm}.

\section{Equations of motion}\label{sec_model}

Consider a satellite model consisting of a rigid body (carrier) with an attached flexible beam. We suppose that the satellite moves on a circular orbit around the Earth and associate
a Cartesian frame $Oxyz$ (body coordinate system -- BCS) with the carrier body.
Another Cartesian frame $Ox'y'z'$ (orbit coordinate system -- OCS) is associated with the orbit such that the axis $Ox'$ is orthogonal to the orbital plane, $Oy'$ is tangent to the orbit, and $Oz'$ is pointing away from the Earth's center.
Let us denote the unit vectors of frames $Oxyz$ and $Ox'y'z'$ by $({\bm e}_1, {\bm e}_2, {\bm e}_3)$ and $({\bm i}, {\bm j}, {\bm k})$, respectively, and introduce the absolute angular velocity vector of the carrier body in BCS as $\bm \omega = \omega_1 {\bm e}_1+ \omega_2 {\bm e}_2 + \omega_3 {\bm e}_3$.

\paragraph{Kinematics.} Let $\omega_0={\rm const}$ be the orbital rate of the considered satellite,
then the absolute angular velocity of OCS is ${\bm \omega}_O = - \omega_0{\bm i}$, and the relative angular velocity of the carrier body with respect to OCS is ${\bm \omega}_r  = \bm \omega + \omega_0{\bm i}$.
The rotation of BCS with respect to OCS is described by the following quaternion equation~\cite[Chap.~8.6]{markley2014fundamentals}:
\begin{equation}\label{eq:QtrnEq}
\begin{aligned}
    & \bm{\dot q}=\frac12q_4{\bm \omega}_r + \frac12 {\bm q} \times {\bm \omega}_r,\\
    & \dot q_4 = - \frac12 \left<{\bm q},{\bm \omega}_r\right>,
    \end{aligned}
\end{equation}
where the dot stands for the time derivative, $\left<\cdot,\cdot\right>$ denotes the inner product in~$\mathbb{R}^3$, ${\bm q}=(q_1,q_2,q_3)^T$ is the vector part of the attitude quaternion and $q_4$ is the scalar part constrained by the condition
$$
q_1^2 + q_2^2 + q_3^2 + q_4^2 = 1.
$$
The unit vectors $(\bm i, \bm j, \bm k)$ of OCS have the following coordinates in BCS~\cite{aerospace22}:
\begin{equation}\label{eq:Vec_k}
{\bm i}=\begin{pmatrix}
q_1^2-q_2^2-q_3^2+q_4^2 \\
2(q_1 q_2-q_3 q_4)\\
2(q_1 q_3+q_2 q_4)
\end{pmatrix},\;
{\bm j}=\begin{pmatrix}
2(q_1 q_2+q_3 q_4) \\
q_4^2-q_1^2+q_2^2-q_3^2\\
2(q_2 q_3-q_1 q_4)
\end{pmatrix},\;
{\bm k}=
    \left(
      \begin{array}{c}
        2(q_1q_3-q_2q_4) \\
        2(q_2q_3+q_1q_4) \\
        q_3^2+q_4^2-q_1^2-q_2^2 \\
      \end{array}
    \right).
\end{equation}


\paragraph{Dynamics of the carrier body.}
It is assumed that the origin of the body coordinate system coincides with the center of mass of the satellite in the undeformed reference state of the beam and that $Ox$, $Oy$, $Oz$ are principal axes of inertia of the satellite considered as a rigid body with the ``frozen'' beam at its reference state.
We follow the restricted formulation of the problem of dynamics (see,~e.g.,~\cite{rumyantsev}) and assume that the overall motion is decoupled into the orbital motion of the center and the spherical motion of the satellite around its center of mass. We also assume that the beam is subjected to small deformations,
so that the mass distribution of the satellite can be described by a constant tensor of inertia, and the rotation of the carrier body is governed by Euler's equations~(cf.~\cite{markley2014fundamentals}):
\begin{equation}\label{eq:EulerEq}
	\hat I\dot{\bm\omega}+\bm\omega\times(\hat I\bm\omega)=\bm u+ {\bm \tau}_g,
\end{equation}
where $\hat I={\rm diag}(I_1, I_2, I_3)$ is the tensor of inertia of the satellite with the ``frozen'' beam, $\dot {\bm \omega} = \dot\omega_1 {\bm e}_1+ \dot\omega_2 {\bm e}_2 + \dot\omega_3 {\bm e}_3$,
$\bm u$ is the control torque,
$
{\bm \tau}_g = 3\,\omega_0^2(\bm k\times \hat I\bm k)
$
is the gravity gradient torque~\cite{wertz}, the components of $\bm k$ in BCS can be computed by formulas~\eqref{eq:Vec_k}.

The three-dimensional control torque $\bm u$ can be implemented, for instance,
in a reaction control system (RCS) with three independent thrusters.
An alternative approach for generating a constrained torque relies on the use of magnetorquers~(see, e.g.,~\cite{ovchinnikov2019survey,silani2005magnetic} and references therein).
If $\bm \mu$ is the magnetic moment generated by magnetorquers, then the control torque acting on the satellite is
$$
{\bm u} = {\bm \mu}\times {\bm B}(t),
$$
where ${\bm B}(t)$ is the Earth magnetic field vector at the corresponding point of the orbit.
The vector function ${\bm B}(t)$ is periodic in the general case, and ${\bm B}(t) = \beta{\bm i} = {\rm const}$ for a satellite in equatorial orbit. The magnetic moment $\bm \mu$ can be regulated by the current applied to a set of electromagnetic coils, so it is possible to treat $\bm \mu\in {\mathbb R}^3$ as the control vector for a satellite model with magnetic actuation.
The method of a fully magnetic attitude control was successfully tested on small satellites having electromagnetic actuation~\cite{sechi2011flight}, and it is considered to be a promising technology based on inexhaustible energy source (using solar panels).

\paragraph{Vibration of the boom.}
We use the Euler--Bernoulli beam model to describe the dynamics of the flexible part (boom) of the satellite.
Let ${\bm w}(\zeta,t)=w_1(\zeta,t){\bm e}_1+w_2(\zeta,t){\bm e}_2+(\zeta+\ell_0){\bm e}_3$ parameterize the center line of the beam in BCS depending on the Lagrangian coordinate $\zeta\in [0,\ell]$ and time $t$.
We assume that the beam of length $\ell$ is attached to the carrier at the point with coordinates $(0,0,\ell_0)$, so that $w_1(\zeta,t)=w_2(\zeta,t)=0$ describe the undeformed reference state of the beam.
The kinetic and potential energy of the beam can be written as
$$K=\frac12\int\limits_0^\ell\ra\left(\dot w_1^2+\dot w_2^2\right)\dxx$$
and
$$U=\frac12\int\limits_0^\ell EI\left((w_1'')^2+(w_2'')^2\right)\dxx,$$
respectively, where the derivatives with respect to $\zeta$ are denoted by a prime,
$\rho$ is density of the beam (mass per unit volume),
$A$ is the cross-sectional area of the beam,
$E$ is the Young's modulus, and $I$ is the cross-sectional moment of inertia of the beam.

Then we construct the Lagrangian $L=K-U$ and apply Hamilton's principle~\cite{berd}: if
$w_i(\zeta,t)$, $i=1,2$ define the motion of the considered mechanical system for
$t\in [t_1,t_2]$, then
\begin{equation}\label{eq:Hamltn-Ostr}
	\delta{\int\limits_{t_1}^{t_2}(K-U)\dd t}+\int\limits_{t_1}^{t_2}\delta{\cal A}\dd t=0
\end{equation}
for every admissible variations $\delta{w}_i(\zeta,t)$ such that $\delta{w}_i(\zeta,t_1)=\delta{w}_i(\zeta,t_2)=0$, $i=1,2$,
where $\int\limits_{t_1}^{t_2}\delta{\cal A}\dd t$ is the work of external forces on $\delta{w}_i$.
The admissible variations $\delta{w}_i(\zeta,t)$ are of class $C^1$ with respect to $t$ and of class $C^2$ with respect to $\zeta$.
We postulate the boundary conditions for the cantilever beam:
\begin{equation}\label{eq:BC}
  \begin{aligned}
    & w_1(0,t)=w_2(0,t)=0,\quad w_1'(0,t)=w_2'(0,t)=0,\\[6pt]
    & w_1''(\ell,t)=w_2''(\ell,t)=0,\quad w_1'''(\ell,t)=w_2'''(\ell,t)=0.
  \end{aligned}
\end{equation}
Performing the integration by parts in~\eqref{eq:Hamltn-Ostr}, we obtain
\begin{equation}\label{eq:var_Lagrangian}
	\delta{\int\limits_{t_1}^{t_2}(K-U)\dd t}=-\int\limits_{t_1}^{t_2}\int\limits_0^\ell\ra\left(\ddot w_1\delta w_1+\ddot w_2\delta w_2\right)\dxx\dd t-\int\limits_{t_1}^{t_2}\int\limits_0^\ell EI\left(w_1^{(4)}\delta w_1+w_2^{(4)}\delta w_2\right)\dxx\dd t.
\end{equation}
Here and further, $w^{(4)}$ stands for the fourth derivative of $w$ with respect to $\zeta$.

In order to describe the motion in the non-inertial frame BCS, we need to introduce the corresponding inertial forces. Let ${\bm F}$ be the inertial force acting on a beam element with the mass $\rho Ad\zeta$. It is well-known that ${\bm F}$ is expressed as the sum of the Euler force~${\bm F}_1$, Coriolis force~${\bm F}_2$, and centrifugal force~${\bm F}_3$~\cite{lanczos}, which can be calculated using the following formulas:
\begin{equation}\label{eq:ext_forces1}
	\tfrac1{\ra\dxx}\,{\bm F}_1 = -\,\bm{\dot\omega}\times{\bm w},\quad
	\tfrac1{\ra\dxx}\,{\bm F}_2 = -2\,{\bm\omega}\times\bm{\dot w},\quad
	\tfrac1{\ra\dxx}\,{\bm F}_3 = -{\bm\omega}\cdot({\bm\omega}\cdot{\bm w})+{\bm w}\cdot{\bm\omega}^2.
\end{equation}
Then the work due to the force ${\bm F}$ is expressed as follows:
\begin{equation*}
	\int\limits_{t_1}^{t_2}\delta{\cal A}\dd t=\int\limits_{t_1}^{t_2}\int\limits_0^\ell\langle{\bm F},\delta{\bm w}\rangle\dxx\dd t,
\end{equation*}
where
\begin{equation}\label{eq:ext_forces2}
    \begin{aligned}
        \tfrac1{\ra\dxx}\,\langle{\bm F},\delta{\bm w}\rangle=
            & \left((\omega_2^2+\omega_3^2)w_1+(\dot\omega_3-\omega_1\omega_2)w_2+2\omega_3\dot w_2-(\dot\omega_2+\omega_1\omega_3)\xx\right)\delta w_1,\\
            & \left((\omega_1^2+\omega_3^2)w_2-(\dot\omega_3-\omega_1\omega_2)w_1-2\omega_3\dot w_1+(\dot\omega_1-\omega_2\omega_3)\xx\right)\delta w_2.
    \end{aligned}
\end{equation}
Calculating the projections of $\int\limits_{t_1}^{t_2}\delta{\cal A}\dd t$ on the $x$ and $y$ axes and applying the fundamental lemma of calculus of variations to~\eqref{eq:Hamltn-Ostr}, we obtain the following differential equations:
\begin{equation}\label{eq:BeamEq}
  \begin{aligned}
    & \ddot w_1+\tfrac{EI}{\ra}w_1^{(4)}=(\omega_2^2+\omega_3^2)w_1+(\dot\omega_3-\omega_1\omega_2)w_2-(\dot\omega_2+\omega_1\omega_3)\xx+2\omega_3\dot w_2,\\[6pt]
    & \ddot w_2+\tfrac{EI}{\ra}w_2^{(4)}=(\omega_1^2+\omega_3^2)w_2-(\dot\omega_3+\omega_1\omega_2)w_1+(\dot\omega_1-\omega_2\omega_3)\xx-2\omega_3\dot w_1,\; \zeta \in [0,\ell].
  \end{aligned}
\end{equation}

\paragraph{Truncated system.}
For further analysis, we consider the case of small oscillations of the beam and small angular velocity of the carrier body, so we will truncate the terms of order higher than 2 with respect to $w_i$, $\omega_j$ in~\eqref{eq:BeamEq}.
Then the resulting system of partial differential equations takes the form
\begin{equation}\label{eq:BeamEq_q}
  \begin{aligned}
    & \ddot w_1+\tfrac{EI}{\ra}w_1^{(4)}=\dot\omega_3 w_2-(\dot\omega_2+\omega_1\omega_3)\xx+2\omega_3\dot w_2,\\[6pt]
    & \ddot w_2+\tfrac{EI}{\ra}w_2^{(4)}=-\dot\omega_3 w_1+(\dot\omega_1-\omega_2\omega_3)\xx-2\omega_3\dot w_1,\quad \zeta \in [0,\ell].
  \end{aligned}
\end{equation}
The derivatives $\dot\omega_j$ can be eliminated from the right-hand side of~\eqref{eq:BeamEq_q} by using the component-wise representation of Euler's equations~\eqref{eq:EulerEq}:
\begin{equation}\label{eq:EulerEqSys}
    \begin{aligned}
        \dot\omega_1 & =\tfrac{I_2-I_3}{I_1}\omega_2\omega_3+ \tfrac{6\omega_0^2(I_2-I_3)}{I_1} (q_1q_4+q_2q_3)(2q_1^2+2q_2^2-1)+\tfrac{u_1}{I_1},\\
        \dot\omega_2 & =\tfrac{I_3-I_1}{I_2}\omega_1\omega_3+\tfrac{6\omega_0^2(I_3-I_1)}{I_2}(q_1q_3-q_2q_4)(2q_1^2+2q_2^2-1)+\tfrac{u_2}{I_2},\\
        \dot\omega_3 & =\tfrac{I_1-I_2}{I_3}\omega_1\omega_2+\tfrac{12 \omega_0^2(I_1-I_2)}{I_3}(q_1 q_4 + q_2 q_3)(q_2q_4-q_1q_3)+\tfrac{u_3}{I_3},
    \end{aligned}
\end{equation}
where we apply the normalization condition $q_4^2 = 1-q_1^2-q_2^2-q_3^2$.



Thus, {\em the nonlinear hybrid system described by ordinary differential equations~\eqref{eq:QtrnEq},~\eqref{eq:EulerEqSys} and
partial differential equations~\eqref{eq:BeamEq_q} with boundary conditions~\eqref{eq:BC} is proposed as a mathematical model of the considered satellite with flexible boom}.

\section{Stabilizing control design}\label{sec_control}
Our goal is to stabilize the undeformed state of the beam $w_1(\zeta,t)=w_1(\zeta,t)=0$ by applying a suitable state feedback control.
 For this purpose, we consider a Lyapunov functional candidate in the form
\begin{equation}\label{eq:E}
        {\cal V} = \tfrac12\intl\left(\dot w_1^2+\dot w_2^2 +\frac{EI}{\rho A}\left((w_1'')^2+(w_2'')^2\right)\right)\dxx.
\end{equation}
By computing the time derivative of $\cal V$ along the trajectories of system~\eqref{eq:BeamEq_q},~\eqref{eq:BC}, we get:
\begin{equation}\label{eq:dot_E_q}
        \dot{\cal V} = (\dot \omega_1 -\omega_2 \omega_3) \cdot \gamma_1(\dot w_2) + (\dot \omega_2 +\omega_1 \omega_3) \cdot \gamma_2(\dot w_1) + \dot \omega_3 \cdot \gamma_3(w_1,w_2,\dot w_1, \dot w_2),
\end{equation}
where the functionals $\gamma_i$ are
\begin{equation}\label{eq:gammas}
    \gamma_1(\dot w_2)=\intl\xx\dot w_2\dxx, \; \gamma_2(\dot w_1)=-\intl\xx\dot w_1\dxx, \; \gamma_3(w_1,w_2,\dot w_1, \dot w_2)=\intl(w_2\dot w_1-w_1\dot w_2)\dxx.
\end{equation}

We will define a feedback control to ensure that the time-derivative $\dot {\cal V}$ along the trajectories of the closed-loop system is negative semidefinite:
\begin{equation}\label{eq:dotE_ctrl}
    \dot{\cal V} = -\nu_1 \gamma_1^2 -\nu_2 \gamma_2^2 -\nu_3 \gamma_3^2 \le 0,
\end{equation}
where $\nu_i>0$ are arbitrary constant gain parameters, $i=\overline{1,3}$.

By exploiting the structure of equations~\eqref{eq:EulerEqSys}, we conclude that formula~\eqref{eq:dotE_ctrl} holds with the control components
\begin{equation}\label{eq:ctrl}
    \begin{aligned}
        & u_1 = -\nu_1 I_1 \gamma_1 + (I_1-I_2+I_3)\omega_2\omega_3 + 6 \omega_0^2 (I_3-I_2) (q_1q_4+q_2q_3)(2q_1^2+2q_2^2-1), \\
        & u_2 = -\nu_2 I_2 \gamma_2 + (I_1-I_2-I_3)\omega_1\omega_3 + 6 \omega_0^2 (I_1-I_3) (q_1q_3-q_2q_4)(2q_1^2+2q_2^2-1), \\
        & u_3 = -\nu_3 I_3 \gamma_3 + (I_2-I_1)\omega_1\omega_2 + 12 \omega_0^2 (I_2-I_1) (q_1 q_4 + q_2 q_3)(q_2q_4-q_1q_3).
    \end{aligned}
\end{equation}

\section{Abstract representation of the closed-loop system}\label{sec_abstract}
To treat the stability problem formally, we will rewrite the obtained closed-loop system in the state space form.
For this purpose, we introduce the real Hilbert space
$$
\begin{aligned}
X = \big\{{\bm\xi} = (w_1,v_1,w_2,v_2,{\bm \omega},{\bm q},q_4)^T\,|& \,w_1,w_2\in H^2(0,\ell),\,v_1,v_2\in L^2(0,\ell),\,{\bm \omega},{\bm q}\in{\mathbb R}^3,\,q_4\in\mathbb R,\\
&\,w_1(0)=w_2(0)=0,\, w_1'(0)=w_2'(0)=0\big\},
\end{aligned}
$$
where the inner product of tuples ${\bm\xi} = (w_1,v_1,w_2,v_2,{\bm \omega},{\bm q},q_4)^T$ and $\bar{\bm\xi} = (\bar w_1,\bar v_1,\bar w_2,\bar v_2,\bar {\bm \omega},\bar{\bm q},\bar q_4)^T$ is defined by
$$
\begin{aligned}
\left<{\bm\xi},\bar{\bm\xi}\right>_X = & \int_0^\ell \Bigl(\rho A (v_1(\zeta)\bar v_1(\zeta)+v_2(\zeta)\bar v_2(\zeta))+EI ( w_1''(\zeta)\bar w_1''(\zeta)+w_2''(\zeta)\bar w_2''(\zeta) )\Bigr)\dxx\\
& + I_1 \omega_1\bar\omega_1 + I_2 \omega_2\bar\omega_2 + I_3 \omega_3\bar\omega_3 + \varkappa( q_1\bar q_1 + q_2\bar q_2 + q_3\bar q_3 + q_4\bar q_4),
\end{aligned}
$$
$\varkappa$ is a positive constant, and $H^2(0,\ell)$ is the Sobolev space of square integrable functions with square integrable first and second derivatives on $[0,\ell]$.

Let $w_1(x,t)$, $w_2(x,t)$, ${\bm \omega}(t)$, ${\bm q}(t)$, $q_4(t)$ satisfy the differential equations~\eqref{eq:QtrnEq},~\eqref{eq:BeamEq_q},~\eqref{eq:EulerEqSys} with boundary conditions~\eqref{eq:BC} and the feedback control~\eqref{eq:ctrl} for $t\in {\cal I}= [0,t_\omega)$, $t_\omega\leqslant +\infty$.
Then a straightforward computation shows that $\bm\xi(t)=(w_1(\cdot,t),\dot w_1(\cdot,t),w_2(\cdot,t),\dot w_2(\cdot,t),\bm\omega(t)$, $\bm q(t),q_4(t))^T\in X$ satisfies the following abstract differential equation:
\begin{equation}\label{eq:abstr}
    \dot{\bm\xi}(t)= F \bm\xi(t),\quad t\in {\cal I},
\end{equation}
where the nonlinear operator $F:D(F)\to X$ is defined by the rule
\begin{equation*}
    F:\bm\xi=\begin{pmatrix}w_1 \\ v_1 \\w_2 \\ v_2 \\ \omega_1 \\ \omega_2 \\ \omega_3 \\ q_1 \\ q_2 \\ q_3 \\ q_4\end{pmatrix} \mapsto F \bm\xi =\begin{pmatrix}
            v_1 \\
            -\tfrac{EI}\ra w_1^{(4)} + 2\omega_3v_2 + {\nu_2}\xx \gamma_2(v_1) - \nu_3 w_2 \gamma_3(w_1,w_2,v_1,v_2) \\
            v_2 \\
            -\tfrac{EI}\ra w_2^{(4)} - 2\omega_3v_1 - \nu_1 \xx \gamma_1(v_2) + \nu_3 w_1 \gamma_3(w_1,w_2,v_1,v_2) \\
            \omega_2\omega_3 - {\nu_1}\gamma_1(v_2) \\
            -\omega_1\omega_3 - {\nu_2}\gamma_2(v_1) \\
            -{\nu_3}\gamma_3(w_1,w_2,v_1,v_2)  \\
            \tfrac12\bigl(\omega_3q_2-\omega_2q_3+(\omega_1+\omega_0)q_4\bigr) \\
            \tfrac12\bigl((\omega_1-\omega_0)q_3-\omega_3q_1+\omega_2q_4\bigr) \\
            \tfrac12\bigl(\omega_2q_1-(\omega_1-\omega_0)q_2+\omega_3q_4\bigr) \\
            -\frac12 \bigl(q_1 (\omega_1+\omega_0) + q_2 \omega_2 + q_3\omega_3\bigr)
          \end{pmatrix},
\end{equation*}
and its domain is
\begin{equation*}
    D(F) = \left\{ \bm\xi\in X:\quad
        \begin{aligned}
            & w_j\in H^4(0,\ell),\, w_j(0)=w_j'(0)=0,\, w_j''(\ell)=w_j'''(\ell)=0, \\
            & v_j\in H^2(0,\ell),\, v_j(0)=v_j'(0)=0,\; j=1,2, \\
            & \omega_i, q_j \in{\mathbb R},\, i=\overline{1,3},\; j=\overline{1,4}
        \end{aligned} \right\}\subset X.
\end{equation*}

In the sequel, {\em equation~\eqref{eq:abstr} will be treated as the abstract formulation of the closed-loop system~\eqref{eq:QtrnEq},~\eqref{eq:BC},~\eqref{eq:BeamEq_q},~\eqref{eq:EulerEqSys},~\eqref{eq:ctrl}}.

\section{Partial stability conditions}\label{sec_thm}

It is easy to see that system~\eqref{eq:abstr} has an equilibrium
\begin{equation}\label{equil}
\hat{\bm \xi}=(0,0,\; 0,0,\; -\omega_0,0,0,\; 0,0,0,1)^T\in X,
\end{equation}
which corresponds to the circular motion of the satellite as an absolutely rigid body around the Earth. In this case, the frames $Oxyz$ and $Ox'y'z'$ coincide, the absolute angular velocity of the satellite is ${\bm \omega}_O = - \omega_0{\bm i}$, and the beam is in its undeformed state.

Let us consider a nonnegative functional $y:X\to {\mathbb R}^+$ such that
\begin{equation}\label{y_fun}
y({\bm \xi}) = \|w_1\|_{H^2(0,\ell)} + \|w_2\|_{H^2(0,\ell)}+\| v_1\|_{L^2(0,\ell)} + \|v_2\|_{L^2(0,\ell)}.
\end{equation}
For the stability analysis, we will make the following assumption.

{\bf Assumption~A1.}
{\em The abstract Cauchy problem~\eqref{eq:abstr} with ${\bm\xi}(0)={\bm \xi}_0\in X$ is well-posed on ${\cal I} = [0,+\infty)$.}

According to the definition of partial stability (cf.~\cite{movchan,rumyantsev,zuyev2003partial}), we call the equilibrium~\eqref{equil} {\em stable with respect to the functional} $y({\bm\xi})$, if:
\begin{itemize}
 \item[(i)] for any
 $\varepsilon>0$, there exists a $\delta(\varepsilon)>0$ such that
 each solution of~\eqref{eq:abstr} satisfying the initial condition
 $\|{\bm\xi}(0) - \hat {\bm\xi}\|_X < \delta(\varepsilon)$ has the property $y({\bm\xi}(t))<
 \varepsilon$ for all $t\ge 0$.
\end{itemize}

We summarize the main result of this paper in the following theorem.

{\bf Theorem~1.}
{\em
Let Assumption~A1 be satisfied for system~\eqref{eq:abstr} with some positive constants $\nu_1$, $\nu_2$, $\nu_3$.
Then the equilibrium ${\bm\xi}=\hat {\bm\xi}$ of~\eqref{eq:abstr} is stable with respect to the functional $y({\bm\xi})$ given by~\eqref{y_fun}.
}

{\em Proof.} Inequality~\eqref{eq:dotE_ctrl} implies that the functional
$$
{\cal V}  ({\bm\xi}) =\tfrac12\intl\left(v_1^2(\zeta)+v_1^2(\zeta)+\frac{EI}{\rho A}\left((w_1''(\zeta))^2+(w_2''(\zeta))^2\right)\right)\dxx.
$$
is nonincreasing along the solutions ${\bm\xi}(t)$ of~\eqref{eq:abstr}.
Moreover, by exploiting Wirtinger's inequality~\cite{hardy1952inequalities}, we conclude that
$$
\alpha_1 y^2 ({\bm\xi}) \le {\cal V} ({\bm\xi}) \le \alpha_2 \|{\bm\xi} - \hat {\bm\xi}\|^2_X\quad \text{for all}\; {\bm\xi}\in X,
$$
with some positive constants $\alpha_1$ and $\alpha_2$.
Therefore, the assertion of Theorem~1 follows from the proof of property~(i) in~\cite[Theorem~2.1]{zuyev2003partial} and~\cite[Theorem~5.2]{movchan}. $\square$

\section{Conclusion}
The feedback control law introduced in~\eqref{eq:ctrl} offers a constructive solution to the considered problem of partial stabilization in the infinite-dimensional state space $X$.
In order to analyze the well-posedness of the Cauchy problem for system~\eqref{eq:abstr}, the theory of $C_0$-semigroups of operators can be applied~\cite{barbu2010nonlinear,krein}.
In this paper, we do not address the issue of well-posedness, deferring its examination to future research endeavors.

\section*{Acknowledgement}
This work was supported by the National Academy of Sciences of Ukraine (budget program KPKVK 6541230, project No. VB-15-18-21/479).

\bibliographystyle{acm}
\bibliography{satellite}

\end{document}